\documentclass[11pt,twoside,a4paper]{article} 
\usepackage[utf8]{inputenc} 
\usepackage{graphicx,color}

\usepackage[bookmarksopen=false,pdftex=true,breaklinks=true,%
      plainpages=false,%
      hyperindex=true,pdfstartview=FitH,%
      pdfpagelabels=true,
      linkcolor=blue,%
      citecolor=blue,urlcolor=red,
      colorlinks=true%
      ]%
   {hyperref}
\usepackage{amssymb}
\usepackage{amsfonts}
\usepackage{amsthm}

\usepackage[all]{xy}
\usepackage[utf8]{inputenc}

\RequirePackage{etoolbox}
\patchcmd{\sectionmark}{\MakeUppercase}{}{}{}
\makeatletter
\renewcommand\@oddhead{Virtual roots\hfill\thepage}
\makeatother
\def\ce#1{\centerline{#1}}
\def\demos{\vskip 1.ex \noindent \sl{Proof. }\rm }
\def \eop {\hbox{}\nobreak\hfill
   \vrule width 1.4mm height 1.4mm depth 0mm
          \par \goodbreak \smallskip}

\def\lcof{{\rm lcof}}

\def\@{\char'76}

\def\eqdf#1{\buildrel{#1}\over =}
\newcommand \C{\mathbb {C}}

\newcommand \Q{\mathbb {Q}} 
\newcommand \R{\mathbb {R}} 
\def\eqdef{\buildrel{\rm def}\over =}

\def\eqdf#1{\buildrel{#1}\over =}
\outer\def\mybeginsection#1\par{\vskip0pt plus.3\vsize
\vskip0pt plus-.3\vsize\bigskip\vskip\parskip  
\message{#1}\noindent {\bf#1}\nobreak\smallskip\noindent}

\def\salce{integral continuous function }
\def\salces{integral continuous functions }
\def\salcs{semialgebraic continuous functions }

\def\boxit#1{\vbox{\hrule\hbox{\vrule\kern3pt
       \vbox{\kern3pt#1\kern3pt}\kern3pt\vrule}\hrule}}

\usepackage{amsthm}

\theoremstyle{plain}
\newtheorem{theorem}{Theorem}[section]
\newtheorem*{theorem*}{Theorem}

\newtheorem{lemma}[theorem]{Lemma}

\newtheorem{proposition}[theorem]{Proposition}

\theoremstyle{definition}

\newtheorem{definition}[theorem]{Definition}

\newtheorem{notation}[theorem]{Notation}

\theoremstyle{remark}

\newtheorem{example}[theorem]{Example}

\begin{document}

\title{Virtual Roots of Real Polynomials} 
\author{Laureano Gonzalez--Vega and Henri Lombardi and Louis Mahé}
\date{2017 (tiny corrections to the original paper, \\note also that the proof of the bound in Theorem \ref{th5.5} is incomplete)}
\maketitle


\begin{abstract}
The fact that
a real univariate polynomial misses some real roots is usually overcame by
considering complex roots, but the price to pay for, is a complete lost of the
sign structure that a set of real roots is endowed with (mutual position on the
line, signs of the derivatives, etc...). In this paper we present  real
substitutes for these missing roots which keep sign properties and which extend
of course the existing roots. Moreover these "virtual roots" are the values of
semialgebraic continuous -- rather uniformly  -- functions defined on the set of
monic polynomials. We present some applications. 
\end{abstract}


\section*{Introduction} 

The problem known as Pierce--Birkhoff Conjecture is the following: take a
real valued continuous function on ${\R}^n$ which is piecewise polynomial,
with a finite number of pieces (i.e. a ``$C_0$--spline"), can you write it down as a finite combination of sup and inf of polynomials? Under this form the problem has been solved for $n\le2$ and the proof for $n=2$ (\cite{Ma}, see also \cite{Del}), uses actually a certain parametrization of the
$1$--dimensional case. Unfortunately this parametrization is not good enough
to get the result for higher dimension and one is still looking for some path in this direction. 

Actually, the proof in the low  dimension case uses the notion
of ``truncation of a polynomial" which is the following: if $u$ is the $r$--th real
zero of the degree $d$ univariate polynomial $P(X)$, the ``$r$--th
truncation" of $P$,
$\phi_{d,r}(P)$ is the function defined as $0$ when $x\le u$ and equal to
$P(x)$ for
$x\ge u$. The essential  point is that $\phi_{d,r}(P)$ is an Inf--Sup definable
function (ISD in short) and that its formal description with sup and inf is
the same
for every other polynomial $Q$ as long as the relative position of the real
roots of
all the successive derivatives of $Q$ is the same as for $P$. This kind of
``local uniformity" makes possible to define
$\phi_{d,r}(P)$ for
multivariate polynomials $P(\underline{X},Y)$, considering $\underline{X}$ as
parameters, as long as $\underline{X}$ belongs to some semi--algebraic set,
precisely described by the sign conditions which define the position of the
zeroes
of the $Y$--derivatives of $P$; and this is sufficient to get the proof in
dimension
$2$. But, for higher dimension, we need more uniformity in the
$1$--dimensional case. In
particular, it would be nice to have this partially defined fonction
$\phi_{d,r}$,
defined everywhere on the parameter space. This is of course impossible in
general: 
the $r$--th real zero alone need not exist for a given value of~$\underline{X}$. Of course, live would be easier if every monic degree $d$
polynomial
would have $d$ real roots!

Actually, the notions of ``virtual root" we are going to introduce
in this paper will
give a good substitute to this unreachable paradise and will, in some sense,
``render hyperbolic every polynomial" (a polynomial is hyperbolic when its
roots are
real). More precisely, we have two classes of ``virtual root functions"
defined on the
set of  degree $d$ monic univariate polynomials of ${\R}[Y]$ (which can be
identified to ${\R}^d$) and one of these classes is the following:

\smallskip \noindent  For every integer $d\ge 1$ and every integer
$0<j\le d$, there is a real valued semi--algebraic continuous function
$\rho_{d,j}$  on
${\R}^d$, such that
$\rho_{d,j}(P)$ is the $j$--th real root of $P$ when $P$ is hyperbolic, and which
satifies in addition the sign conditions we expect for an actual $j$--th
root. For
example $\rho_{d,j}(P)\le\rho_{d-1,j}(P')\le\rho_{d,j+1}(P)$ if $P'$ is the
derivative of $P$.

\smallskip \noindent Then, once we have our hands on the $r$--th virtual
root
$\rho_{d,r}(P)(\underline{X})$ of a degree $d\ge r$ monic polynomial
$P(\underline{X},Y)$ everywhere on the parameter space, the next step towards a
solution of Pierce--Birkhoff Conjecture would be to construct the  ``$r$--th
virtual truncation" of $P$ as an ISD function coinciding with $P$ for
$Y\ge\rho_{d,r}$ and
``going to zero as fast as possible" for $Y\le\rho_{d,r}$,
and giving of course the actual truncation in case $\rho_{d,r}$ is an actual
root. This is not yet completely worked out and should appear in a next future.
Nevertheless, as early applications of these notions, we prove here the two
following results:
\begin{itemize}
\item [1)] a continuous version of Thom's lemma,
\item [2)] the closure under Sup and Inf of the ring  generated by the virtual
roots is the integral closure of the polynomial ring ${\R}[X_1,\ldots ,X_n]$
inside the real valued continuous functions on ${\R}^n$.
\end{itemize}

\noindent The paper is organized in the following way.
\renewcommand\contentsname{\relax}
\vspace{-3em}
\tableofcontents
%

\bigskip
In what follows, we have choosen to work over the real numbers $\R$, but
everything can be worked out over any real closed field. Even more, if the
polynomials we start with have their coefficients in a subfield $\bf K$ of
a real closed field $\bf R$, every new constructed polynomial has also its
coefficients in this field $\bf K$.

\section{General tools} \label{sec1}

\noindent As we said in the abstract, we want to define on the set of monic
univariate real polynomials,  some collections of ``virtual root" functions,
extending everywhere the actual root functions in such a way that some sign
conditions are preserved. There are essentially two ways to distinguish a
given real root of a polynomial out of the others: one is the rank of this
root, the other is the collection of the signs taken at this root by the
derivatives.

The main idea is the simple following observation: suppose $P$ is a
parametrized polynomial in one variable and we are following some
particular real root along the parameters. If for some value of the
parameters this root disappears, then it becomes a root of the derivative and
this becomes our ``virtual root". But in both cases, actual or virtual root,
the root realizes the local minimum of the absolute value of $P$, and this is
the key observation.

So, we are going to
consider two sets of such root functions, called respectively ``$r$--th 
virtual root" and ``Thom's virtual roots". In the first case we want to
preserve the rank of a given root among the others. In the second case we
try to preserve the sign that every derivative of $P$ takes on a given root,
but it is a bit more complicated.

The main tool to define these functions is the following one.
\begin{definition} \label{defi1.1}
 We identify the set of monic degree $d$ polynomials of ${\bf
R}[X]$  to ${\R}^d$, and~$P$ will be understood as a polynomial or as a
point in ${\R}^d$ as well. Let ${\cal S}_d$ be the closed
$\Q$--semialgebraic set defined by: $${\cal S}_d=\{ (a,b,P):a\leq
b,\;\;\deg(P)=d,\;\;\forall x,y\in [a,b]\quad P'(x)P'(y) \geq 0\}$$ and
${\cal R}_d$ be the semialgebraic function defined on ${\cal S}_d$  by:
$${\cal R}_d(a,b,P)=z\quad\hbox{\rm such that}\quad \vert P(z)\vert=
\min\{\vert P(u)\vert:u\in [a,b]\}$$
\end{definition}

An easy verification shows that the function ${\cal R}_d$ satisfies the
following  equality:
$${\cal R}_d(a,b,P)=\cases{\hfill a = b \hfill& if a = b\cr
                           \hfill a \hfill& if $(P(b)-P(a))\cdot
P(a)\geq 0$\cr
                           \hfill b \hfill& if $(P(b)-P(a))\cdot
P(b)\leq 0$\cr
                           \hbox{the real root of $P$ in $(a,b)$}
\hfill& otherwise}\eqno(\star)$$

\begin{proposition} \label{prop1.2}
If $d$ is a nonnegative integer then we have.
\begin{enumerate}
\item If  $(a,b,P)\in{\cal S}_d$ and  $P$ has a real root $z$ 
on   $[a,b]$ then ${\cal R}_d(a,b,P)=z$,
\item the function ${\cal R}_d$ is continuous on ${\cal S}_d$,
\item if $(a,b,P)\in{\cal S}_d$ then the number  
$x={\cal R}_d(a,b,P)$ can be characterized by the following
inequalities.
$$\matrix{a\leq x\leq b\hfill&&\cr\cr
          (x-a)P(a)(P(b)-P(a))\leq 0&&(x-a)P(x)(P(b)-P(a))\leq
0\cr\cr
          (b-x)P(b)(P(b)-P(a))\geq 0&&(b-x)P(x)(P(b)-P(a))\geq
0\cr}$$
\end{enumerate}
 
\end{proposition}

\demos

\noindent Items \emph{1} and \emph{3} are easy considering the different cases
appearing  in the formula ($\star$). Next we prove item \emph{2} which is no more 
than proving that the real root of a monotone polynomial in an interval varies
continuously with the coefficients. Let
$(a,b,P)$ be an element in ${\cal S}_d$ and $\epsilon$ a strictly
positive element of $\R$. We search for a $\delta$ giving the
continuity of the function ${\cal R}_d$. 

If $b-a\leq \epsilon/2$ then taking $\delta=\epsilon/2$ we have:
$$\!\vert a-a'\vert+\vert b-b'\vert+\vert
P-R\vert<\delta\,\Longrightarrow
\,
\vert x-x'\vert\leq \max\{b,b'\} - \min\{a,a'\} \leq
\vert b-a\vert+\vert a-a'\vert+\vert b-b'\vert <
\epsilon$$ 
with $(a',b',R)\in{\cal S}_d$, $x={\cal
R}_d(a,b,P)$ and 
$x'={\cal R}_d(a',b',R)$.

If $b-a\geq\epsilon/2$ and $x={\cal R}_d(a,b,P)$ then,
writing
$\alpha$ for
$+1$ or $-1$ according to the sign of  $P(b)-P(a)$, we
consider  three cases.

\smallskip
\begin{itemize}
\item [{$\bullet$}] If $x<a+\epsilon/2$ then  $\alpha\cdot
P(a+\epsilon/2)>0$. For  a sufficiently small variation $\delta$ of
$(a,b,P)$ in
${\cal S}_d$, the real number $\alpha\cdot P(a+\epsilon /2)$
remains strictly positive, the variation of $a$ is smaller than
$\epsilon /2$ and
${\cal R}_d(a,b,P)$ remains on the interval
$[a,a+\epsilon /2)$.
\item [{$\bullet$}] If  $x > b-\epsilon / 2$,  we proceed  the
same  way as in the previous case. 
\item [{$\bullet$}] If  $a+\epsilon /2 \leq x \leq b-\epsilon / 2$,
then 
$\alpha\cdot P(x-\epsilon / 4)<0<\alpha\cdot P(x+\epsilon / 4)$.
For  a sufficiently small variation $\delta$ of $(a,b,P)$ in
${\cal S}_d$, 
$\alpha\cdot P(x-\epsilon /4)$ remains $<0$, $\alpha\cdot
P(x+\epsilon /4)$ remains strictly positive, and the variations of
$a$ and $b$ are smaller than
$\epsilon /4$. So, ${\cal R}_d(a,b,P)$ remains in the open
interval 
$(x-\epsilon /4 , x+\epsilon / 4)$. 
\end{itemize}
\vspace{-2em} \eop

\medskip
Next we generalize the definition of ${\cal R}_d$ to the cases
$a=-\infty$ or
$b=+\infty$. This is achieved by considering the semialgebraic sets:
$${\cal S}_{d,+}=\{(a,P):\forall x\in[a,+\infty)\quad\quad P'(x)\geq
0\}$$
$${\cal S}_{d,-}=\{(b,P):\forall x\in(-\infty,b]\quad\quad
(-1)^{d-1}P'(x)\geq  0\}$$ defining on ${\cal S}_{d,+}$:
$${\cal R}_d(a,+\infty,P)={\cal R}_d(a,\max(a,1+
\sup_{i=0,\ldots,d-1}\{\vert a_i\vert\}),P)$$  and defining on
${\cal S}_{d,-}$: $${\cal R}_d(-\infty,b,P)= {\cal R}_d(\min(b,-1-
\sup_{i=0,\ldots,d-1}\{\vert a_i\vert\}),b,P)$$
\medskip

\begin{notation} \label{nota1.3}
If $Q$ is a univariate polynomial then $Q^{(i)}$ will
denote the 
$i$--th  derivative of $Q$, with $Q^{(0)}=Q$, $\deg(Q)$ will be the
degree of
$Q$ and $\lcof(Q)$ its leading coefficient. In order to be able to use the
identification between $\R^d$ and the set of monic degree $d$ polynomials,
we define: $$Q^{[i]}={{Q^{(d-i)}}\over{\lcof(Q^{(d-i)})}}$$ as the
normalized derivative of $Q$ of degree $i$. We define  $Q^\star$ as the
product of all normalized derivatives $Q^{[i]}$ of $Q$ ($Q$ included).\end{notation}

\section{The $r$--th virtual root} 

\noindent Let $P\in{\R}[X]$ be a monic degree $d$ polynomial. For every integer
$r$ such that $0<r\le d$, we want to define a function $\rho_{d,r}$ on ${\R}^d$
having the following properties.
\begin{enumerate}
\item $\rho_{d,r}$ is a continuous semialgebraic function on ${\R}^d$.
\item If $P$ is hyperbolic and $u\in {\R}$ is the $r$--th real root of the polynomial $P$, then $u=\rho_{d,r}(P)$.
\item $\rho_{d,r}(P)\le\rho_{d-1,r}(P'/d)\le\rho_{d,r+1}(P)$.
\end{enumerate}

The restriction to monic polynomials is not really essential: we could
be satisfied with polynomials such that the leading coefficient never
vanishes, but then we would loose some uniformity in the continuity of
$\rho_{d,r}$ (see section {\bf 5}). But without loss of generality, we
may as well replace monic by ``quasi--monic", meaning that the leading
coefficient is~$\geq 1$. Anyway, for simplicity, we will do
everything with monic polynomials.

\begin{definition} \label{defi2.1}
Let $P(x)=x^d-(a_{d-1}x^{d-1}+\ldots+a_0)$ be a monic
polynomial in $\R[x]$. For $d\ge 0$ and for any integer $j$, we define 
$\rho_{d,j}(P)$ in the following inductive way.
\begin{enumerate}
\item [{$\bullet$}] If $j\le 0$, we put $\rho_{d,j}(P)=-\infty$,
\item [{$\bullet$}] if $j>d$, we put $\rho_{d,j}(P)=\infty$,
\item [{$\bullet$}] if  $d>0$ and $1\le j\le d$, we define
$$ \rho_{d,j}(P)\eqdf{\rm def}
 {\cal R}_d(\rho_{d-1,j-1}(P'/d),\rho_{d-1,j}(P'/d),P)$$
\end{enumerate}

In particular if $P=X-a$ then $\rho_{1,1}(P)=a$. Let's also define the sets:
$${\rm U}_{d,j}(P)\eqdf{\rm def} \{ \alpha \in {\R}:
\rho_{d,j-1}(P) <\alpha < \rho_{d,j}(P)\}$$  
\end{definition}
These sets are open
intervals when they are not empty, and they are empty in particular for
$j<0$ and $j>d$. For simplicity, we will often write $U_{d,j}(P')$ and
$\rho_{d,j}(P')$ instead of the corresponding terms with $P'/d$.
Applying Proposition \ref{prop1.2}, it is easy to prove by induction the following
proposition:
\begin{proposition} \label{prop2.2}
 
For $d>0$ and $0<j\le d$, the functions
$\rho_{d,j}$ are \salces on ${\R}^d$ defined over  $\Q$ and they are
roots of  the  polynomial  $P^*$. On the other hand every root of 
$P$  is equal to some 
$\rho_{d,j}(P)$.
\end{proposition}

Let us quote here the basic properties of these $\rho_{d,r}$ which make $P$
looking like hyperbolic with respect to the virtual roots:

\begin{proposition} \label{prop2.3}
For $d>0$, the functions $\rho_{d,r}$ have the following
properties.
\begin{enumerate}
\item $\forall r \,\rho_{d,r}(P)\le\rho_{d-1,r}(P')\le\rho_{d,r+1}(P)$.
\item Every monic degree $d$ polynomial has $d$ virtual roots (possibly
equal).
\item $(-1)^{d+r}P(x)>0$ for $x\in U_{d,r+1}(P)$.
\end{enumerate} 
\end{proposition}

\demos

\noindent Items \emph{1} and \emph{2} are just from the definition. For \emph{3}, we make an
induction on $d$. Anyway, there is something to prove only when the interval
$(\rho_{d,r}(P),\rho_{d,r+1}(P))$ is not empty, so we may assume $0\le r\le
d$.

If $d=1$, it is easily checked. If $d>1$, we have
$$\rho_{d-1,r-1}(P')\le\rho_{d,r}(P)\le\rho_{d-1,r}(P')\le\rho_{d,r+1}(P)\le
\rho_{d-1,r+1}(P')$$
we consider two cases:
\begin{enumerate}
\item[{$\star$}] if $\rho_{d,r}(P)=\rho_{d-1,r}(P')$, then
$U_{d,r+1}(P)\subseteq U_{d-1,r+1}(P')$ and we know by hypothesis that
$(-1)^{d+r}P'<0$ on $U_{d-1,r+1}(P')$. So  $(-1)^{d+r}P$ is
decreasing on $U_{d,r+1}(P)$. As $\rho_{d,r+1}(P)$ realizes the minimum
of $|P|$ on $U_{d-1,r+1}(P')$, we get that  $(-1)^{d+r}P>0$ on
$U_{d,r+1}(P)$.

\item[{$\star$}] If $\rho_{d,r}(P)<\rho_{d-1,r}(P')$, then $(-1)^{d-1+r-1}P'>0$ on
$U_{d-1,r}(P')$ and $(-1)^{d+r}P$ is increasing on this interval. As 
$\rho_{d,r}(P)$ realizes the minimum of $|P|$ on this interval,
$(-1)^{d+r}P$ must be positive on $U_{d-1,r}(P')\cap
U_{d,r+1}(P)\not=\emptyset$, and must be so on the whole of $U_{d,r+1}(P)$,
for it cannot change sign on this interval.
\end{enumerate}
\vspace{-2em}\eop

\section{Thom's virtual roots}  

 Let $P(x)=x^d-(a_{d-1}x^{d-1}+\ldots+a_0)$ be a  monic polynomial in 
${\R}[x]$. Thom's lemma says in particular that if we fix the sign (in
the large sense) of every derivative of $P$, we get a set containing at most
one root of $P$. The virtual roots we are going to build up are real numbers
$x_\sigma$ indexed by a list of signs
$\sigma=[\sigma_0,\ldots,\sigma_{d-1}]$ having the property that when  the
$d-1$ non trivial derivatives of $P$ take the sign given by the list
$\sigma$ at some real root of $P$, then $x_\sigma$ is precisely this root.
Of course, it may happen that the sign condition on the derivatives produce
an empty set: in that case the point $x_\sigma$ cannot satisfy the sign
conditions (although it might be an actual root of $P$).

\begin{notation} \label{nota3.1}
We shall denote by $\sigma=[\sigma_0,\ldots,\sigma_d]$ a
list of  signs, $\sigma_i\in\{+,-\}$. The ``length"
${\rm lg}(\sigma)$   will be $d$ and  $\sigma_0$ 
will  always be equal   to $+$. The convenience of this $\sigma_0$ will
appear later. Concerning the list $\sigma$ we introduce the following
symbols, for $i=1\ldots d$: 
$$\widetilde{\sigma}_i=\cases{>& if
$\sigma_i=+$\cr <& if $\sigma_i=-$\cr} \quad
\overline{\sigma}_i=\cases{\geq& if $\sigma_i=+$\cr \leq& if
$\sigma_i=-$\cr}
\quad
\sigma^{(i)}=[\sigma_0,\ldots,\sigma_{d-i}]
\quad
\sigma^{[i]}=[\sigma_0,\ldots,\sigma_{i}]$$
The basic semialgebraic open set
$$\{\alpha\in{\R}:P^{[1]}(\alpha)\;\widetilde{\sigma}_1\;
0,\ldots, P^{[d]}(\alpha)\;\widetilde{\sigma}_d\; 0\}$$ will be
denoted by ${\rm U}_{\sigma}(P)$ and the basic semialgebraic closed
set
$$\{\alpha\in{\R}:P^{[1]}(\alpha)\;\overline{\sigma}_1\; 0,\ldots,
P^{[d]}(\alpha)\;\overline{\sigma}_d\; 0\}$$ by ${\rm
F}_{\sigma}(P)$.
\end{notation}

With the previous notations Thom's Lemma (see \cite{BCR} for a proof) can
be stated  in the following terms.

\begin{theorem}[Thom's Lemma] \label{th3.2}
If the closed set ${\rm F}_{\sigma}(P)$ is not  empty
then it is a  closed interval or a point, and its interior is always
${\rm U}_{\sigma}(P)$. Moreover every finite  endpoint of the
interval ${\rm F}_{\sigma}(P)$ is a root of some $P^{(j)}$. 
\end{theorem}
\begin{definition} \label{defi3.3}
Suppose $\deg P={\rm lg}(\sigma)=d$ and $\epsilon\in\{+,-\}$.
Here  we assume that ${\rm
F}_{\sigma}(P)$ is not empty. The two endpoints of  ${\rm F}_{\sigma}(P)$ will
be denoted by: $\tau^{\epsilon}_{\sigma}(P)$ with
$\epsilon=+$ for the right endpoint and $\epsilon=-$ for the left endpoint. 
\end{definition}

There are two special cases where the interval ${\rm F}_{\sigma}(P)$ is never
empty and one of its endpoints is infinity: $$\sigma=[+,+,+,\ldots,+],\quad
\epsilon=+\quad\Longrightarrow \quad \tau^{\epsilon}_{\sigma}(P)=+\infty$$
$$\sigma=[+,-,+,-,+,\ldots],\quad \epsilon=-\quad\Longrightarrow
\quad \tau^{\epsilon}_{\sigma}(P)=-\infty$$ Excepting the two
infinity cases, the symbols 
$\tau^{\epsilon}_{\sigma}$  represent semialgebraic functions
partially  defined on ${\R}^d$.
In the following two cases, the symbol $\tau^{\epsilon}_{\sigma}$
provides a semialgebraic function defined  on the whole ${\R}^d$:
$$\sigma=[+,+,+,\ldots,+],\quad \epsilon=-\quad\Longrightarrow
\quad \tau^{\epsilon}_{\sigma}(P)=\max\{\alpha\in{\R}: 
P^{*}(\alpha)=0\}$$
$$\sigma=[+,-,+,-,+,\ldots],\quad \epsilon=+\quad\Longrightarrow
\quad\tau^{\epsilon}_{\sigma}(P)=\min\{\alpha\in{\bf
R}:P^{*}(\alpha)=0\}$$

Let us introduce the function, also partially defined on ${\R}^d$, denoted
by $\rho_{\sigma^{(1)}}(P)$ and called {\sl actual Thom's root} 
which is
 defined as the only real root of $P$ inside the closed interval 
$[\tau_{[\sigma_0,\ldots,\sigma_{d-1}]}^-(P'/d),
\tau_{[\sigma_0,\ldots,\sigma_{d-1}]}^+(P'/d)]$ when the endpoints
of the interval are defined (possibly equal)
 and when such a root exists.   The function $\rho_{\sigma^{(1)}}$, when
defined, verifies the following  equality:
$$\rho_{[\sigma_0,\ldots,\sigma_{d-1}]}(P)=
  \tau_{[\sigma_0,\ldots,\sigma_{d-1},\sigma_{d-1}]}^-(P)=
  \tau_{[\sigma_0,\ldots,\sigma_{d-1},-\sigma_{d-1}]}^+(P)$$ It is
clear from the definition that every root of $P$ in ${\R}$, can be
represented by some of these symbols. The functions $\tau$ and $\rho$ will be
extended as \salcs to the  whole of ${\R}^d$ in an inductive way.

\begin{definition} \label{defi3.4}
If the degree $d$ of $P$ is equal to $1$ we define:

\smallskip
\ce{$\rho_{[+]}(x-a)\eqdef \tau_{[+,-]}^+(x-a)\eqdef
\tau_{[+,+]}^-(x-a) \eqdef  a$}
\smallskip
\ce{$\tau_{[+,-]}^-(x-a)\eqdef -\infty\qquad \tau_{[+,+]}^+(x-a)\eqdef
+\infty$}
\smallskip
\noindent If all the functions $\rho$ and $\tau$ for degree $d-1$
are known then the  definitions for degree $d$ are:
\smallskip
\ce{$
  \rho_{[\sigma_0,\ldots,\sigma_{d-1}]}(P)\eqdf{\rm def}
  \tau_{[\sigma_0,\ldots,\sigma_{d-1},\sigma_{d-1}]}^-(P)\eqdf{\rm
def}
  \tau_{[\sigma_0,\ldots,\sigma_{d-1},-\sigma_{d-1}]}^+(P)\eqdf{\rm
def}$}
\smallskip
\ce{$\eqdf{\rm def}{\cal
R}_d(\tau_{[\sigma_0,\ldots,\sigma_{d-1}]}^-(P'/d),
  \tau_{[\sigma_0,\ldots,\sigma_{d-1}]}^+(P'/d),P)$}
\smallskip
\ce{$\tau_{[\sigma_0,\ldots,\sigma_{d-1},\sigma_{d-1}]}^+(P)\eqdf{\rm
def}
  \tau_{[\sigma_0,\ldots,\sigma_{d-1}]}^+(P'/d)
\qquad\qquad
\tau_{[\sigma_0,\ldots,\sigma_{d-1},-\sigma_{d-1}]}^-(P)\eqdf{\rm
def}
  \tau_{[\sigma_0,\ldots,\sigma_{d-1}]}^-(P'/d)$} 
\end{definition}

Remark that if $\epsilon\cdot\sigma _{d}= +$ then
$\tau_{\sigma}^{\epsilon}(P)=\tau_{\sigma ^{(1)}}^{\epsilon}(P'/d)$ and if
$\epsilon\cdot\sigma _{d}= -$ then $\tau_{\sigma}^{\epsilon}(P)=\rho_{\sigma
^{(1)}}(P)$. So we see inductively that excepting the infinity cases each
function $P \mapsto \tau_{\sigma}^{\epsilon}(P)$ is equal to some function
$P \mapsto \rho_{\sigma ^{[j]}}(P^{[j+1]})$, where $j<d$ depends only on
$\sigma $ and $\epsilon$. We then get the following proposition:

\begin{proposition} \label{prop3.5}~
\begin{enumerate}
\item The above defined functions $\rho_{\sigma}$ and
$\tau_{\sigma}^{\epsilon}$ are defined on the whole of ${\R}^d$ and are
extensions of the partial functions introduced in Definition \ref{defi3.3}.
\item The functions $\rho_{\sigma}$ are integral,
$\Q$--semialgebraic and continuous on ${\R}^d$,  and  verify,
for every monic polynomial $P$
 of degree $d$, the equality $P^*(\rho_{\sigma}(P))=0$.
\end{enumerate}
 
\end{proposition}

\demos

\noindent The proof is easy  by induction on the degree, using Proposition
\ref{prop1.2} for the continuity in~1), and that $\rho_\sigma=\rho_\sigma (P)$ for
the generic monic degree $d$ polynomial to show it is integral in~2).\eop

In order to understand these functions $\rho_\sigma$, it is convenient to
introduce the following definition:

\begin{definition} \label{defi3.6}
For every monic polynomial $P$ and every $\sigma$ of length $d$, we
define $G_{\sigma}(P)=[\tau_{\sigma}^-(P),\tau_{\sigma}^+(P)]$. By
construction, this closed interval (maybe a point) depends continuously on
$P$ and coincide with $F_\sigma$ when the latter is not empty. 
\end{definition}

The main properties of $G_{\sigma}$ are summarized below:

\begin{proposition} \label{prop3.7}
 The interval $G_\sigma$ has the following properties.
\begin{enumerate}
\item[{a-.}] ${\rm G}_{[+,+]}(x-a)=[a,+\infty]$ and ${\rm
G}_{[+,-]}(x-a)=[-\infty,a]$,
\item[{b-.}] the interval ${\rm
G}_{[\sigma_0,\ldots,\sigma_{d-1}]}(P')$  is the union of the two
intervals 
$${\rm G}_{[\sigma_0,\ldots,\sigma_{d-1},-\sigma_{d-1}]}(P)\quad
\hbox{and} \quad {\rm
G}_{[\sigma_0,\ldots,\sigma_{d-1},\sigma_{d-1}]}(P)$$  with the
right endpoint of the first one equal to the left endpoint of the 
second one, 
\item[{c-.}] if  ${\rm G}_{\sigma}(P)$
is not reduced  to a point then ${\rm G}_{\sigma}(P)=F_\sigma(P)$.
\end{enumerate}
 
\end{proposition}

\demos

\noindent Point \emph{b-} comes right from Definition \ref{defi3.4}. For \emph{c-}, it is clear by
definition and case examination, that if ${\rm G}_{\sigma}(P)$
is not reduced  to a point then ${\rm G}_{\sigma}(P)$ is the set of
points $\alpha$ in
${\rm G}_{\sigma^{(1)}}(P')$ such that $P(\alpha)\overline{\sigma}_d
0$. But then by induction on $d$, we may assume ${\rm
G}_{\sigma^{(1)}}(P')={\rm F}_{\sigma^{(1)}}(P')$ and so
$G_\sigma(P)=F_\sigma(P)$.
\eop

We can now  understand better the functions $\rho_\sigma$ themselves:  the 
virtual Thom's roots are actual Thom's roots of some derivative:

\begin{proposition} \label{prop3.8}
For every monic polynomial $P$ (resp. $Q$) of degree $d$
(resp. $d+1$) and~$\sigma$ of length $d-1$ (resp. $d$), each $\rho_{\sigma}(P)$
(resp. $\tau^\epsilon_\sigma(Q)$) is equal to an actual Thom's root
$\rho_{\sigma^{[r-1]}}(P^{[r]})$ for some $r$. 
\end{proposition}

\demos

\noindent By definition of $\tau$, it is sufficient to do it for $\rho_\sigma$.
If $d=1$, it is the definition. If $d>1$, let $u=\rho_{\sigma}(P)$. By
construction, $u\in G_\sigma(P^{[d-1]})$ and if $u$ is not an endpoint of this
interval, it is a root of $P$. But in that case, by Proposition \ref{prop3.7}
c-, it is the  actual Thom's root of~$P$ coded by $\sigma$. So we may
assume  $u$ is  an endpoint of this interval. Let $r$ be the smallest integer
such that $u$ is an endpoint of $G_{\sigma^{[r+1]}}(P^{[r]})$. By Proposition \ref{prop3.7} b-, $u$ is inside $G_{\sigma^{[r]}}(P^{[r-1]})$, and by the same argument as
above, must be a root of $P^{[r]}$, coded by $\sigma^{[r]}$. 
\eop

In the case of $r$--th virtual roots the general pattern is
quite easy: there is at most~$d$ virtual roots in degree $d$,
naturally ordered and there is  generically exactly $d$
distincts such virtual roots (realized in specializing to hyperbolic
polynomials). On the contrary, the situation for virtual Thom's roots
is not so clear: How many generic such roots do we have and how are
they mutually ordered? Is there some specialization that gives
the $2^{d-1}$ a priori possible $\rho_\sigma$? We have the two 
following propositions. 

\begin{proposition} \label{prop3.9}
For every $d\ge 1$ and every $\sigma$ of length $d-1$, there is
a  real polynomial~$P$ of degree $d$  such that $\rho_\sigma(P)$ is an
actual Thom's root of $P$. 
\end{proposition}

\demos

\noindent It is sufficient to show that for any  sign condition $\sigma$ of
length $d-1$ there is a  real polynomial~$P$ of degree $d$ having a
root inside $U_\sigma(P)$. In degree $1$ there is nothing to do and \hbox{if
$d>1$}, by induction we may assume that there exists $Q$ of degree
$d-1$ having a root in $U_{\sigma^{(1)}}(Q)$, making
$U_\sigma(P)\not=\emptyset$ for any antiderivative $P$ of $Q$. Adjusting
the constant term, it is then easy to find such a $P$ having a
root in $U_\sigma(P)$.\eop

\medskip

Of course, this implies that there are $2^{d-1}$ distinct generic
Thom's virtual roots.

\begin{proposition} \label{prop3.10}
Let $s(d)=1+d(d-1)/2$.
\begin{enumerate}
\item[{a)}] Every monic degree $d$ polynomial has at most $s(d)$ distinct  Thom's virtual roots.
\item[{b)}]  For every $d\ge 1$ there exists a monic polynomial $P$ which has 
$s(d)$ distinct Thom's virtual roots. 
\end{enumerate} 
\end{proposition}

\demos

\noindent By definition of the $\rho_\sigma(P)$ (length($\sigma)=d-1$), there is
exactly one $\rho_\sigma$ in each $G_\sigma(P')$, and in particular
there is also exactly one in each non empty $F_\sigma$. But the
non empty $F_\sigma$ make a partition of $\R$ and their endpoints
are zeroes of $P'^\star$: the number of intervals is then bounded by
one more than the number of roots of $P'^\star$, which gives a).

Now we show that this bound $s(d)$ for  non empty
$F_\sigma(P)$  is effectively obtained. Choose~$P$ such that $P^\star$
is hyperbolic without multiple roots. If $d=1$ or $2$, it is clear. If
$d>2$, assume the number of intervals for $P'$
(determined by sign conditions on $P^{(i)}$, $i\ge 2$) is $s(d-1)$,
then there are  $d-1$ intervals actually cutted in two by the
$d-1$ roots of $P'$ and the number of non empty $F_\sigma$ is exactly
$s(d-1)+d-1=s(d)$. Let us show that the~$\rho_\sigma$ corresponding
to these $s(d)$ intervals produce $s(d)$ different real numbers: if
two such~$\rho_\sigma$ would be equal, they would correspond to a
common end of two consecutive intervals, realizing the minimal of
the absolute value of $P$ on the union of these two intervals. We
have two cases.
\begin{enumerate}
\item[{1)}] $|P|$ has a positive minimum at that point, but then cannot
be hyperbolic,
\item[{2)}] The point is a root of $P$, but is also a root of $P'^\star$
as an end of a $F_\sigma$ , giving a double root to $P^\star$:
contradiction.
\end{enumerate}
\vspace{-2em} \eop

\medskip

Informations about how the functions $\rho_\sigma$ are ordered are
summarized in the next proposition.

\begin{proposition} \label{prop3.11}
Assume that $\deg(P)=d$ and $\sigma$ is a list of signs ($+$ or
$-$) with ${\rm lg}(\sigma)=d-1$.
\begin{enumerate}
\item [a-.] If  $\mu$ is a list with length $k-1$ (with $d\geq k\geq 1$)
different from $\sigma$, then the comparison, by 
$\geq$ or  $\leq$,  between $ \rho_{\sigma}(P)$ and $\rho_{\mu}(P^{[k]})$ is
given by the  following rule involving only the signs in $\sigma$ and $\mu$.
\\
 If $i$ is the first index such 
that $\sigma_i\neq \mu_i$ (if  $\mu$ is an initial segment of $\sigma$ then
$i={\rm lg}(\mu) +1$) then  the sign of $\rho_{\sigma}(P)-\rho_{\mu}(P^{[k]})$
is equal to $\sigma_{i-1}\cdot \sigma_i$.
\item [b-.] If $u$ is an element of ${\R}$ then the comparison between $u$ and
$\rho_{\sigma}(P)$ is given by ``the same" rule than in (a) using the sign of
$P^{[j]}(u)$ instead $\mu_j$.
\\
 If $i$ is the first index such 
that $\sigma_i\neq {\rm sign}(P^{[i]}(u))$   and $i<d$ then the sign of
$\rho_{\sigma}(P)-u$ is equal to $\sigma_{i-1}\cdot \sigma_i$. If $\sigma_i
= {\rm sign}(P^{[i]}(u))$ for $i=1,\ldots,d-1$  then the sign of
$\rho_{\sigma}(P)-u$ is equal to $-\sigma_{i-1}\cdot {\rm sign}(P(u))$
%
%
\end{enumerate}
\end{proposition}

\demos

\noindent Part (a) is a direct consequence of the formal
construction of the symbols $\rho$. Part (b) comes from (a) when
$u=\rho_{\sigma}(P)$. If
$u\neq\rho$ then the result is clear when the considered symbol
$\rho_\sigma$ corresponds to an actual Thom's root of $P$ coded   by 
$\sigma$. As any $\rho_\sigma$  is an actual
Thom's root of some derivative $P^{[j+1]}$ coded   by 
$\sigma^{[j]}$  we compare $u$ with
$\rho_{\sigma^{[j]}}(P^{[j+1]})$ as in the previous case  (details left to
the reader).\eop

\newpage

\section{Examples}

\begin{example} \label{exa4.1}
We present the picture of the complete situation of $r$--th virtual roots
$\rho_{4,j}(x):=\rho_{4,j}(P)(x)$ of the polynomial
$P(x,y)=((x-1)^2+(y+1)^2-2)((x+1)^2+(y-1)^2-2)$ considered as a
polynomial in $y$ parametrized by $x$. In the picture we can see
 the union of two circles corresponding to the zeroes of $P$,
 a cubic corresponding to the zeroes of $P'_y$,
 an ellipsis corresponding to the zeroes of $P''_y$
and  the $y$--axis being the zero locus of $P^{(3)}_y$.
The number $j$ on the picture denotes $\rho_{4,j}(x)$ and
$\rho_{4,2}(x)$ has been drawn in thick.
\vspace{-10.5cm}
\begin{center}
\includegraphics*[width=13.5cm]{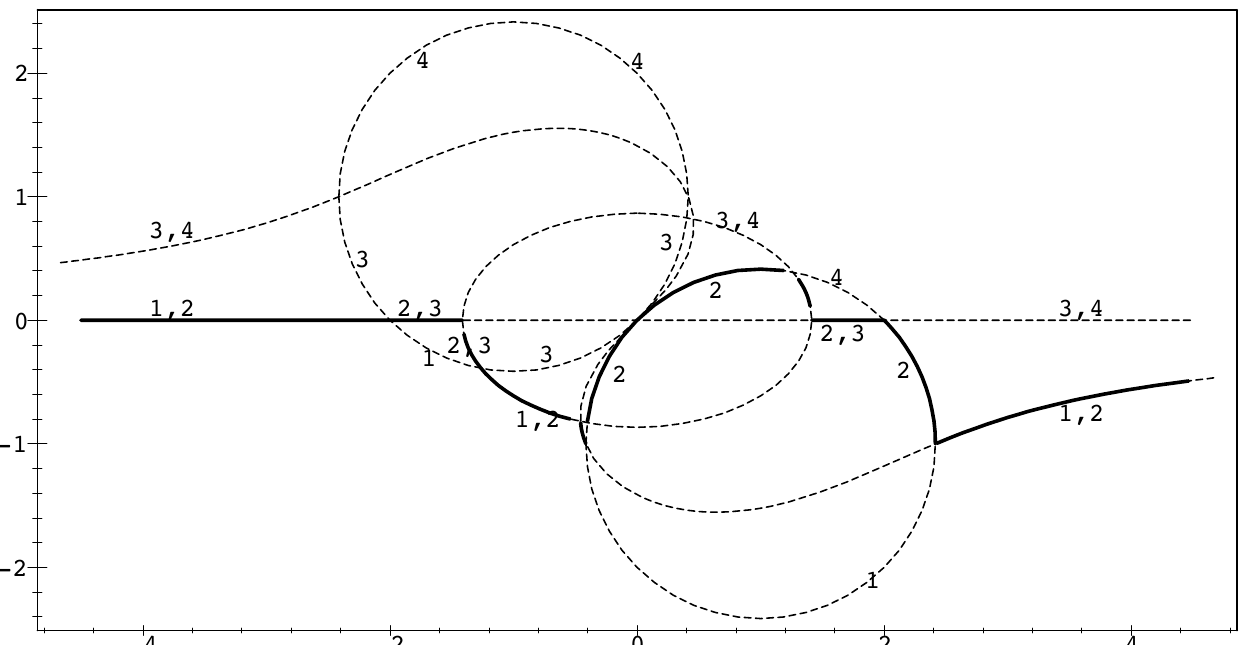}
\end{center}

\end{example}

\begin{example} \label{exa4.2} 
The table below gives the complete situation of the $\rho_\sigma$
until the degree $5$, and is easy to extend to any degree:

\bigskip
\ce{\boxit{\vbox{\offinterlineskip
\halign{&#\cr 
$\;\downarrow$ degree of the polynomial derivative\cr}
\smallskip\hrule
\halign{&#\vrule &\hfil#\hfil&\hfil#\hfil&\hfil#\hfil&\hfil#\hfil&\hfil#\hfil&\hfil
#\hfil&\hfil#\hfil
&\hfil#\hfil&\hfil#\hfil&\hfil#\hfil&\hfil#\hfil&\hfil#\hfil&\hfil#\hfil&\hfil#\hfil
&\hfil#\hfil&\hfil#\hfil&\hfil#\hfil&\hfil#\hfil&\hfil#\hfil&\hfil#\hfil&\hfil#\hfil
&\hfil#\hfil&\hfil#\hfil&\hfil#\hfil&\hfil#\hfil&\hfil#\hfil&\hfil#\hfil&\hfil#\hfil
&\hfil#\hfil&\hfil#\hfil&\hfil#\hfil&\hfil#\hfil\vrule\cr\strut
$\;0\;$&&&&&&&&&&&&&&&&$\;+\;$&&&&&&&&&&&&&&&\cr\strut
$\;1\;$&&&&&&&&$\;-\;$&&&&&&&&$\bullet$&&&&&&&&$\;+\;$&&&&&&&\cr\strut
$\;2\;$&&&&$\;+\;$&&&&$\bullet$&&&&$\;-\;$&&&&\vrule&&&&$\;-\;$&&&&$\bullet$&&&&$\;+\;$
&&&\cr\strut$\;3\;$&&$\;-\;$&&$\bullet$&&$\;+\;$&&\vrule&&$\;+\;$&&$\bullet$&&$\;-\;$&&
\vrule&&$\;+\;$&&$\bullet$&&$\;-\;$&&\vrule&&$\;-\;$&&$\bullet$&&$\;+\;$&\cr\strut$\;4
\;$&$\;+\;$&$\bullet$&$\;-\;$&\vrule&$\;-\;$&$\bullet$&$\;+\;$
&\vrule&$\;-\;$&$\bullet$&$\;+\;$&\vrule&$\;+\;$&$\bullet$&$\;-\;$&\vrule
&$\;-\;$&$\bullet$&$\;+\;$&\vrule&$\;+\;$&$\bullet$&$\;-\;$&\vrule&$\;+\;$&$\bullet$&$
\;-\;$&\vrule&$\;-\;$&$\bullet$&$\;+\;$\cr\strut$\;5\;$      
&$\bullet$&\vrule&$\bullet$&\vrule&$\bullet$&\vrule&$\bullet$&\vrule&$\bullet$&\vrule&
$\bullet$&\vrule
&$\bullet$&\vrule&$\bullet$&\vrule&$\bullet$&\vrule&$\bullet$&\vrule&$\bullet$&\vrule&
$\bullet$&\vrule&$\bullet$&\vrule&$\bullet$&\vrule&$\bullet$&\vrule&$\bullet$\cr}}}}
\medskip

Every point $\bullet$ in the table denotes a function
$\rho_{\sigma}$ where  the list $\sigma$ is obtained reading from
the top until the considered point 
$\bullet$. If we want to add a line to the table, we do it in such a way that
each sign of the bottom line subdivides in two, the first sign of
the two being the opposite of the existing sign. 
In the previous table it is easy to find
some evident  incompatibilities:

\smallskip \noindent [{$\bullet$}] In degree $3$ it is impossible to have the
symbols
$\rho_{[+,-,-]}$ and $\rho_{[+,+,-]}$ representing the real roots of a
 polynomial because we would have a polynomial with two
consecutive simple roots giving the same sign to the derivative,

\smallskip \noindent [$\bullet$] In degree $4$ we get two incompatibilities with
the same type than the previous one, $\rho_{[+,-,+,+]}$ with
$\rho_{[+,-,-,+]}$ and 
$\rho_{[+,+,-,-]}$ with $\rho_{[+,+,+,-]}$.

\smallskip \noindent [{$\bullet$}] Again in degree $4$ a stronger new type of
incompatibility appears: it is impossible to have simultaneously non
empty the two consecutive intervals ${\rm F}_{[+,-,-,-]}(P)$ and ${\rm
F}_{[+,+,-,+]}(P)$. If 
${\rm F}_{[+,-,-,-]}(P)$ and ${\rm F}_{[+,+,-,+]}(P)$ were non empty
then the polynomial $P'$  would decrease from $-$ to $+$.

\smallskip \noindent [{$\bullet$}] If, for example,  $\rho_{[+,-,-,-]}(P)$ is an
actual Thom's root,  then the
interval ${\rm G}_{[+,+,-,+]}(P)$ is formed by only one  point.
Moreover, in this case, the roots coded by $[+,+,-,-]$ and $[+,+,+,-]$
can  not exist for $P$.

\noindent An exhaustive analysis of the previous table allows to find, by
similar  arguments, all the possible simultaneous Thom's codings for
the real roots of the same polynomial. 

\end{example}
\begin{example} \label{exa4.3} 
 Max and Min are $r$--th
root functions and Thom's root functions:
\vspace{-.5em}
$$\max\{a_1,\ldots,a_k\}=\rho_{k,k}(\prod\nolimits_{i=1}^k(x-a_i))=\rho_{[+,+,+,\ldots,+]}
(\prod\nolimits_{i=1}^k(x-a_i))$$
\vspace{-1.5em}
$$\min\{a_1,\ldots,a_k\}=\rho_{k,1}(\prod\nolimits_{i=1}^k(x-a_i))=
\rho_{[+,-,+,-,\ldots]}(\prod\nolimits_{i=1}^k(x-a_i))$$ 
\noindent The $n$--th root function can be described as:
\vspace{-.5em}
$$\root n \of
{\max(a,0)}=\rho_{n,n}(x^n-a)=\rho_{[+,+,+,\ldots,+]}(x^n-a)$$ 
\end{example}

\begin{example} \label{exa4.4} 
{\it Root functions for a polynomial of degree $3$}

\noindent We consider the polynomial $P=x^3+3px+2q$. The complement of
$pq(p^3+q^2)=0$ in the plane $(p,q)$ has six connected components,
$\{A_i:1\leq i\leq 6\}$, obtained by giving strict signs to $p$, $q$ and
$p^3+q^2$. The border of these open sets will not be considered
because the  root functions extend there continuously.

\vspace{-.5cm}
\begin{center}
\includegraphics*[width=7.5cm]{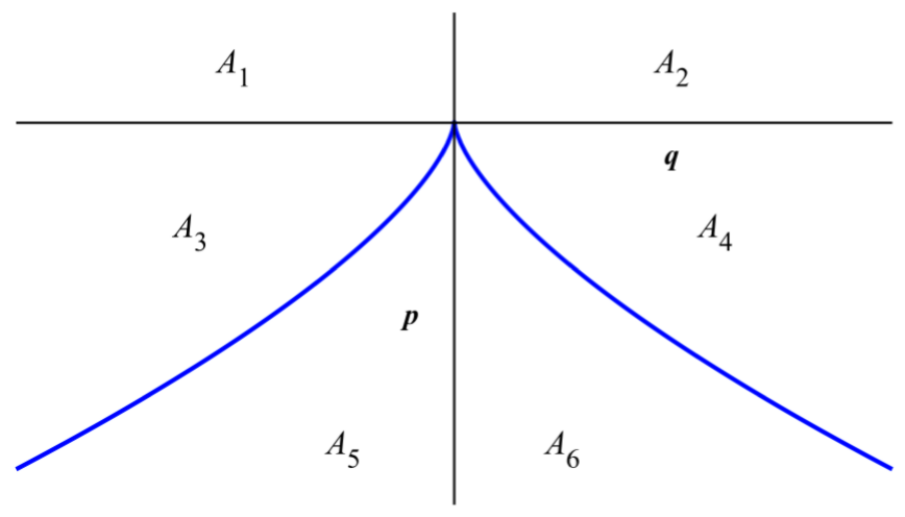}
\end{center}
\vspace{-.5cm}

Inside every $A_i$ each of the four Thom's root functions has a 
fixed expression as an actual Thom's root of $P$ or one of
its derivatives. This fact is shown in the following table:
$$
\matrix{\rho_{[+,+,+]}(P)=\left\{\matrix{
\hbox{$\rho_{[+,+,+]}(P)=$ the biggest real root of $P$ if 
$(p,q)\in A_1\cup A_3\cup A_5\cup A_6$}\hfill\cr
\hbox{$\rho_{[+]}(P^{[1]})=0$ if $(p,q)\in A_2$}\hfill\cr
\hbox{$\rho_{[+,+]}(P^{[2]})=$ the positive real root of $P'$ (i.e.
$\sqrt{-p}$)  if $(p,q)\in A_4$}\hfill\cr}\right.\hfill\cr\cr\cr
\rho_{[+,-,+]}(P)=\left\{\matrix{
\hbox{$\rho_{[+,-,+]}(P)$= the smallest real root of $P$ if 
$(p,q)\in A_2\cup A_4\cup A_5\cup A_6$}\hfill\cr
\hbox{$\rho_{[]}(P^{[1]})=0$ if $(p,q)\in A_1$}\hfill\cr
\hbox{$\rho_{[+,-]}(P^{[2]})=$ the negative real root of $P'$ (i.e.
$-\sqrt{-p}$)  if $(p,q)\in A_3$}\hfill\cr}\right.\hfill\cr\cr\cr
\rho_{[+,+,-]}(P)=\left\{\matrix{
\hbox{$\rho_{[+,+,-]}(P)=$ the intermediate real root of $P$ if 
$(p,q)\in A_6$}\hfill\cr
\hbox{$\rho_{[+]}(P^{[1]})=0$ if $(p,q)\in A_1\cup A_2\cup A_3\cup
A_5$}\hfill\cr \hbox{$\rho_{[+,+]}(P^{[2]})=$ the positive real root of
$P'$ (i.e. $\sqrt{-p}$)  if $(p,q)\in
A_4$}\hfill\cr}\right.\hfill\cr\cr\cr \rho_{[+,-,-]}(P)=\left\{\matrix{
\hbox{$\rho_{[+,-,-]}(P)=$ the intermediate real root of $P$ if 
$(p,q)\in A_5$}\hfill\cr
\hbox{$\rho_{[+]}(P^{[1]})=0$ if $(p,q)\in A_1\cup A_2\cup A_4\cup
A_6$}\hfill\cr \hbox{$\rho_{[+,-]}(P^{[2]})=$ the negative real root of
$P'$ (i.e. $-\sqrt{-p}$)  if $(p,q)\in A_3$}\hfill\cr}\right.\hfill\cr}$$
\end{example}

\newpage
\section{Links and common properties} 

In this section we are going to examine the relationship between the
two kinds of virtual roots, and the properties they share. A
question that comes first in mind is the following: is it possible to
express one set of virtual roots in term of the other? Proposition
\ref{prop5.2} below shows that the $\rho_\sigma$ can be expressed in terms of
$\rho_{d,j}$, but the converse is not yet known. Let us start with a
definition.

\begin{definition} \label{defi5.1}
Let  $\sigma$ be a list of length $d>0$, always with $\sigma_0=+$.
 We define $j(\sigma)$ as $1+\sum_{i=0}^d
(1+\sigma_i\sigma_{i+1})/2$ (the number of ``no sign change" in 
$\sigma$ plus one). For
instance $j([+,-,-])=2$ or $j([+,+,-,-,-])=4$. 
 
\end{definition}

Then we have the following result.

\begin{proposition} \label{prop5.2}
Let $P$ be a degree $d$ monic polynomial and $\sigma$ a list of
length $d$. Then,
\begin{enumerate}
\item If ${\rm U}_{\sigma}(P)\neq\emptyset$ then 
${\rm U}_{\sigma}(P)\subseteq {\rm U}_{d,j(\sigma)}(P)$
\item If  $\rho_{\sigma}(P)$ is an actual Thom's root, then
$\rho_{\sigma}(P)=\rho_{d,j(\sigma)}(P)$
\item In general, for  $\rho_{\sigma}(P)$, we have the
following expression, which allows to express inductively the
$\rho_{\sigma}$ functions as sup--inf combinations of the
$\rho_{d,j}$ functions:
$$\rho_{\sigma}(P)=\max\{\tau_{\sigma}^-(P'),\min\{\tau_{\sigma}^+(P'),
\rho_{d,j(\sigma)}(P)\}\}$$
\end{enumerate} 
\end{proposition}

\demos 

\noindent Let us prove 1) by induction on $d$. If $d=1$, everything is easy. If
$d>1$, if $U_\sigma(P)\not=\emptyset$, it is also the case for
$U_{\sigma^{(1)}}(P')$ and so
$U_{\sigma^{(1)}}(P')\subseteq U_{d-1,j(\sigma^{(1)})}$ by induction.
By Proposition~\ref{prop3.7}, we have two cases to consider:

\smallskip
\noindent {a)} $U_{\sigma^{(1)}}(P')=U_\sigma(P)$ and in that case there is no
zero of $P$ inside $U_\sigma(P)$ and in particular
$\rho_{d,j(\sigma)}$ and $\rho_{d,j(\sigma)-1}$ are outside
$U_\sigma(P)$ (otherwise they would be in
$U_{d-1,j(\sigma^{(1)})}(P')$ and they would be zeroes of $P$ inside
$U_\sigma(P)$). So ${\rm U}_{\sigma}(P)\subseteq {\rm
U}_{d,j(\sigma)}(P)$.

\smallskip
\noindent {b)} $F_{\sigma^{(1)}}(P')$ is the union of two non empty intervals
$F_{\sigma^{(1)},-\sigma_{d-1}}(P)\cup
F_{\sigma^{(1)},\sigma_{d-1}}(P)$, meaning the common endpoint is a
zero of $P$ inside $U_{d-1,j(\sigma^{(1)})}(P')$: it must be
$\rho_{d,j(\sigma^{(1)})}(P)$. So the left interval
$U_{\sigma^{(1)},-\sigma_{d-1}}(P)$ is contained in
$U_{d,j(\sigma^{(1)})}(P)=U_{d,j(\sigma^{(1)},-\sigma_{d-1})}(P)$ and
the right one $U_{\sigma^{(1)},\sigma_{d-1}}(P)$ is contained in
$U_{d,1+j(\sigma^{(1)})}(P)=U_{d,j(\sigma^{(1)},\sigma_{d-1})}(P)$,
which proves 1).\smallskip

It is not hard to see that 1) implies 2): if  $\rho_{\sigma}(P)$ is an
actual Thom's root, then $U_{\sigma^{(1)}}(P')$ is not empty and is
contained in $U_{d-1,j(\sigma^{(1)})}(P')$. The same is true for the
closed corresponding intervals and the only zero of $P$ in
$F_{d-1,j(\sigma^{(1)})}(P')$ is then
$\rho_{d,j(\sigma)}(P)=\rho_{\sigma}(P)$.

It is now easy to show 3): if $U_{\sigma^{(1)}}(P')=\emptyset$, then
$G_{\sigma^{(1)}}(P')$ is a point and the formula in 3) for
$\rho_\sigma(P)$ gives that point. If it is not empty, the same
arguments as above show that, if $\rho_\sigma(P)$ is inside this
open set, it must be equal to $\rho_{d,j(\sigma)}(P)$, and if it is
an endpoint of $F_{\sigma^{(1)}}(P')$, then one of the intervals 
$U_{\sigma^{(1)},-\sigma_{d-1}}(P)$ or
$U_{\sigma^{(1)},\sigma_{d-1}}(P)$ is not empty and so contained in the
corresponding $U_{d,j(\sigma)}(P)$, $\sigma$ being one of the two
possible extensions of $\sigma^{(1)}$. But then the formula gives the
end point of $F_{\sigma^{(1)}}(P')$ corresponding to
$\rho_\sigma(P)$.
Finally, use the remark following Definition \ref{defi3.4} in order to replace in the
formula \ref{prop5.2}~\emph{3} $\tau_{\sigma}^{\epsilon}(P')$ by some $\rho_{\sigma
^{[j]}}(P^{[j+1]})$
(where $ j< d-1$ depends only on $\sigma $ and $\epsilon$).
\eop

\smallskip 
We have already proved that the virtual roots are continuous in the
coefficients of the polynomials, but we know a little more: on a
given compact ball of $\R^d$, they are of course uniformly continuous
and we can compute the modulus of continuity in terms of the radius
of the ball. Let us start with a notation:

\begin{notation} \label{nota5.3}
Let  $X$  be a complete metric space.  We shall denote by
${\rm Ms_k}(X)$ the metric space of multisets with  $k$ elements in 
$X$, i.e. the complete metric space obtained from  $X^k$ with the
semidistance
$${\rm d}_{\rm Ms}((x_i)_{i=1,\ldots ,k},(y_i)_{i=1,\ldots ,k})
\eqdef
\min_{\lambda \in {\rm S}(k)}\{\max_{i=1,\ldots ,k }
d(x_i,y_{\lambda (i)})\}.$$
\end{notation}

Then we need the following lemma:

\begin{lemma} \label{lem5.4}
Let  $U$  be a convex set in  ${\bf R}^n$, $X$ a
metric space, 
$f$: $U \longrightarrow X$   a  continuous function  and  
$F:\;U\longrightarrow\;{\rm Ms_k}(X)$   a uniformly continuous
function with modulus of uniform continuity
$\omega(\epsilon)$. Assume that for all $u \in U$, $f(u) \in  F(u)$.
Then  $f$ admits as  modulus of uniform continuity the function:
$\epsilon \longmapsto \omega(\epsilon /2k).$ 
\end{lemma}

\demos

\noindent We assume w.l.o.g. that $U$ is the unit interval and
$u=0$. We start with  $\epsilon$ and search for  $\delta$ such that
for all $u'\in (0,
\delta)$, we have  $d(f(u),f(u')) < \epsilon$. Let
$\epsilon'=\epsilon /k$,
$\delta = \omega(\epsilon /2k)$ and $u'\in (0, \delta)$. 

Either all $F(u)$ is in the open ball 
${\rm B}_X(f(u),(k-1)\epsilon')$, and then 
$d(f(u),f(u')) < (k-1)\epsilon'+\epsilon /2k < \epsilon$. Or there exists
$j<k-1$ such that $F(u)$ is contained in the disjoint union of the open ball
 ${\rm B}_X(f(u),j\epsilon')$ and of the complement  $X-{\rm
B}_X(f(u),(j+1)\epsilon')$.  Then, for all $t \in [0,u']$, the set
$F(t)$ is contained in the disjoint union of the open ball 
 ${\rm B}_X(f(u),j\epsilon'+\epsilon /2k)$ and of the complement of
the corresponding closed ball 
 $X-{\rm \overline{B}}_X(f(u),(j+1)\epsilon'-\epsilon /2k)$.   So,
by connexity, the point  $f(t)$ must remain in the first ot these
two  disjoint open sets and  $d(f(u),f(u')) < j\epsilon'+\epsilon
/2k < \epsilon$.
\eop
\medskip

Then we have

\begin{theorem}[Root functions local uniform continuity] \label{th5.5}
When the $\vert a_i\vert^{d-i}$ are  bounded by
$M\geq 1$, a modulus of uniform continuity $\omega(M,\epsilon)$ 
(i.e. a function giving $\delta$ from $\epsilon$ in the definition
of uniform continuity, with the $l_1$ norm in $\R^d$) for the
functions 
$\rho_{\sigma}(a_{d-1},\ldots,a_{0})$ and  
$\rho_{d,j}(a_{d-1},\ldots,a_{0})$  is
$$\omega(M,\epsilon)=2M\bigg({{\epsilon}\over{d(d+1)(2d-1)M}}\bigg)^d$$ 
\end{theorem}

\demos

\noindent  A modulus of
uniform continuity $\omega (M,\epsilon )$ for the functions
$\rho_{\sigma}$ and $\rho_{d,j}$ is obtained using the following 
technical result appearing in  \cite[appendix A, page 276]{Ost}.
\smallskip

{\narrower\smallskip\noindent\sl The multiset root function for monic degree
$d$ polynomials:
$$\matrix{\C^d&\longrightarrow&{\rm Ms}_{d}(\C)&\cr
          P&\longmapsto&\{\alpha\in\C:P(\alpha)=0\}\cr}$$
(considering multiplicity) when the $\vert a_i\vert^{d-i}$ are 
bounded by $M\geq 1$, admits the following modulus of local uniform
continuity, with the $l_1$ norm in $\C^d$:
$$2M\bigg({{\epsilon}\over{2M(2d-1)}}\bigg)^d$$\smallskip}

Applying the lemma with the Ostrowski modulus for the
multiset union of complex roots of the polynomials $P$ and its
derivatives (the modulus of uniform continuity for the multiset of
zeros of $P$ is also good for its derivatives), we get the
theorem.  \eop
\medskip

Remark that Ostrowski's bound and Lemma \ref{lem5.4} allow to determine
explicitely  a modulus of local uniform continuity for any integral
semialgebraic  continuous function by merely regarding the vanishing
monic polynomial for   the considered function.

\section{Applications} 
In this section, we conclude with two applications. The first one is
the following continuous version of Thom's lemma:

\begin{theorem}[A continuous version for Thom's Lemma] \label{th6.1}
Let $d$ be an integer $\geq 1$ and 
$\sigma=[\sigma_1,\ldots,\sigma_d]$ a list of elements in $\{+,-\}$. 
We shall consider the monic polynomials with degree $d$ as points of 
${\R}^d$. If we define the sets of ${\R}^d$:
$${\rm W}_{\sigma}=\{P\in {\R}^d: {\rm
F}_{\sigma}(P)\neq\emptyset\}\qquad
  {\rm V}_{\sigma}=\{P\in {\R}^d: {\rm
U}_{\sigma}(P)\neq\emptyset\}$$ then the following statements are
verified:
\begin{enumerate}
\item  ${\rm W}_{\sigma}$ is a connected and closed
$\Q$--semialgebraic  set whose interior is  ${\rm V}_{\sigma}$.
\item  ${\rm V}_{\sigma}$ is a connected and open
$\Q$--semialgebraic  set whose closure is  ${\rm W}_{\sigma}$.
\item  For every $P$ in ${\rm W}_{\sigma}$ the set ${\rm
F}_{\sigma}(P)$  is a non--empty closed interval and every finite
end--point of ${\rm F}_{\sigma}(P)$ is an
\salce  of $P$ and a root of  $P^*$. 
\item  Only two cases with an infinity end--point can appear:
$$\sigma=[+,+,\ldots,+]\longrightarrow +\infty\qquad
  \sigma=[+,-,+,-,+,\ldots]\longrightarrow -\infty$$
\end{enumerate} 
\end{theorem}

\demos

\noindent  Parts (3) and (4) are clear after the detailed study on
the sets ${\rm F}_{\sigma}(P)$ made in the previous sections. The following
equivalences:
$${\rm F}_{\sigma}(P)\neq\emptyset\quad\Longleftrightarrow
\quad\tau^+_{\sigma}(P)\in {\rm F}_{\sigma}(P)\quad\qquad {\rm
U}_{\sigma}(P)\neq\emptyset\quad\Longleftrightarrow\quad
{{\tau^+_{\sigma}(P)+\tau^-_{\sigma}(P)}\over 2}\in {\rm
U}_{\sigma}(P)$$ allow us to show that ${\rm W}_{\sigma}$ is a closed
$\Q$--semialgebraic set and that ${\rm V}_{\sigma}$ is an open
$\Q$--semialgebraic set. 

Now we suppose  w.l.o.g.  that $\sigma_d=\sigma_{d-1}=+$ and that we
are not  in an infinity case. For a degree $d-1$ polynomial $R$ we
define:
$$R_1(x)=d\int_0^x R(t)dt\qquad \psi(R)=R_1(\tau^+_{\sigma^{(1)}}(R))$$ A
simple verification provides the following description for the sets ${\rm
W}_{\sigma}$ and ${\rm V}_{\sigma}$:
$${\rm W}_{\sigma}=\{P:{\rm F}_{\sigma^{(1)}}(P')\neq\emptyset,\;\;
\psi (P'/d)\geq -P(0)\}={\rm W}_{\sigma^{(1)}}\times{\bf
R}\cap\{P:\psi\circ\pi(P)\geq-P(0)\}$$
$${\rm V}_{\sigma}=\{P:{\rm U}_{\sigma^{(1)}}(P')\neq\emptyset,\;\;
\psi (P'/d)> -P(0)\}={\rm V}_{\sigma^{(1)}}\times{\bf
R}\cap\{P:\psi\circ\pi(P)>-P(0)\}$$ where $\pi$ is the projection:
$$\matrix{\pi\colon&{\R}^d&\longrightarrow&{\R}^{d-1}\cr
&P&\longmapsto& P'/d\cr}$$ Proceeding by induction on $d$ we obtain
the remaining claims in (1) and (2) because ${\rm W}_{\sigma}$ and 
${\rm V}_{\sigma}$ are cylinders bounded from below by the continuous
semialgebraic function $\psi\circ\pi+P(0)$ and whose base is a
semialgebraic set verifying the conditions in (1) and (2) by
induction hypothesis.\eop

\medskip

The second question we want to address here is the following:
what kind of functions do we get if we take the closure under inf--sup of the
set of functions $\rho_\sigma$ or $\rho_{d,j}$? If we take a given continuous
function on $\R^n$ which is integral over the $n$ variable polynomials $\R
[X_1,\ldots,X_n]$, it is annihilated by a monic polynomial
$Q(\underline{X},Y)$ in $Y$ with coefficients in $\R[X_1,\ldots,X_n]$, and
piecewise on $\R^n$, it is a precise real root of $Q(\underline{X},Y)$ (in
terms of $r$--th roots or Thom's roots), but in general, it does not admit a
global description as Inf--Sup of the virtual roots of $Q(\underline{X}, Y)$. A
very simple example is the following:

\begin{example} \label{exa6.2} 
Take $Q(X,Y)=Y^2-X^2$, and $f(X)=X$. If we had a description of $f$ as
Inf--Sup of virtual roots of $Q$, it would depend only on $X^2$, and so would
be the same for $X>0$ and $X<0$.  Of course, we have other nice descriptions for
$f$!. But it means that if we want do describe \salces as Inf--Sup of virtual
roots, we have to use other polynomials that $Q$. That is Theorem \ref{th6.4}.
\end{example}

\begin{definition} \label{defi6.3}
If $\rho$ is either a $r$--th root or Thom's
root function on $\R^d$, we define functions on $\R^n$ in filling each
occurrence of $\rho$ with a polynomial in $n$ variables. Let us call
``polyroots in $n$ variables" these functions on $\R^n$ (in both cases), and
``Inf--Sup of polyroots" the functions obtained in taking finite Infima and
Suprema of such functions. 
\end{definition}

Then we get the following:

\begin{theorem} \label{th6.4}
The closure  of polyroots in $h$
variables under sum, Inf and Sup (in both cases of polyroots) is the 
integral closure of $\R [X_1,\ldots,X_h]$ in the ring of continuous functions
on $\R^h$. 
\end{theorem}

\demos

\noindent It is clear that the Inf--Sup of sums of polyroots in $h$ variables
are continuous and integral over the polynomial ring $\R [X_1,\ldots,X_h]$, so
the only thing to prove is the converse.
 Let $f:{\R}^h\longrightarrow{\R}$ be an \salce and
$Q(x_1,\ldots,x_h,y)$ a  polynomial in ${\R}[x_1$,$\ldots$,$x_h$,$y]$,
$y$--monic, with  degree $d$ in $y$, and verifying:
$$Q(\alpha_1,\ldots,\alpha_h,f(\alpha_1,\ldots,\alpha_h))=0
\qquad\quad \forall (\alpha_1,\ldots,\alpha_h)\in{\R}^h $$ We
shall denote $\underline{x}=(x_1,\ldots,x_h)$ and write:
$$Q(\underline{x},y)=y^d+\sum\nolimits_{k=0}^{d-1} Q_k(\underline{x})y^k$$

Let $g_1,\ldots,g_m$ be the virtual root functions corresponding to degree
$d$ ($m=d$ in case of $\rho_{d,j}$ and $m=2^{d-1}$ in case of Thom's
roots).  Then for every $i\in\{1,\ldots,m\}$ the function defined by:
$$l_i(\underline{x})=g_i(Q_{d-1}(\underline{x}),\ldots, Q_0(\underline{x}))$$
is a polyroot. After  these definitions it is clear that
the function $$\prod\nolimits_{i=1}^m (f(\underline{x})-l_i(\underline{x}))$$ is zero
everywhere.

Next, for every $i\in\{1,\ldots,m\}$, we introduce the closed 
semialgebraic set: 
$$F_i=\{(\alpha_1,\ldots,\alpha_h)\in{\R}^h:
f(\alpha_1,\ldots,\alpha_h)=l_i(\alpha_1,\ldots,\alpha_h)\}$$ whose
interior will be denoted by $U_i$.

Applying the Finiteness Theorem we describe every $U_i$ as a finite
union  of basic semialgebraic open sets , i.e. by strict sign
conditions over  polynomials in ${\R}[x_1,\ldots,x_h]$. Let
$\{P_j:j\in J\}$ be the  family of polynomials appearing in such
description and $\{P_j:j\in K\}$ the  family  obtained  completing
the previous one until obtaining a separating family.

Finally we  consider the non--empty  open sets obtained in giving
strict signs to the polynomials in $\{P_j:j\in K\}$. This family will
be denoted by $\{V_n:n\in N\}$.  As our  family of  polynomials is
separating then the closed semialgebraic set obtained replacing  in
the description for $V_n$ the conditions $<$ by $\leq$ and  the
conditions $>$ by $\geq$ is the closure of $V_n$. Moreover, after the
definition of the $V_n$'s it is clear that they are disjoint:
$n\neq p\quad \Longleftrightarrow\quad V_n\cap V_p=\emptyset$

The conclusion of the theorem will be obtained in constructing  a sum of 
Inf--Sup of  polyroots  equal to $f$ over the  union
of the sets $V_n$ (which is dense in ${\R}^h$).

For every $n\in N$ let $i_n$ be such that $V_n\subseteq U_{i_n}$:
this  implies that the function $f$, over $V_n$, is equal to
$l_{i_n}$. Now  we construct for every pair $(n,p)$ with $n\neq p$
an Inf--Sup of 
polyroots  $v_{n,p}$  verifying the following
conditions:
$$\forall\underline{\alpha}\in V_n \quad v_{n,p}(\underline{\alpha})
\geq f(\underline{\alpha})=l_{i_n}(\underline{\alpha})$$
$$\forall\underline{\alpha}\in V_p \quad v_{n,p}(\underline{\alpha})
\leq f(\underline{\alpha})=l_{i_p}(\underline{\alpha})$$ 
If $i_n=i_p$  we define $v_{n,p}=l_{i_n}$. So, without loss of
generality  we can assume that $(n,p)=(1,2)$,  $f=l_1$ on $V_1$ and
$f=l_2$ on $V_2$. Let
$W_1$ and $W_2$ be the  closures of $V_1$ and $V_2$ and write
(w.l.o.g.):
$$V_1=\{\underline{\alpha}\in{\R}^h:
P_1(\underline{\alpha})>0,\ldots,P_r(\underline{\alpha})>0,
\ldots,P_s(\underline{\alpha})>0\}$$
$$V_2=\{\underline{\alpha}\in{\R}^h:
P_1(\underline{\alpha})>0,\ldots,P_r(\underline{\alpha})>0,
P_{r+1}(\underline{\alpha})<0,\ldots,P_s(\underline{\alpha})<0\}$$
This allows to derive the following descriptions for $W_1$ and $W_2$:
$$W_1=\{\underline{\alpha}\in{\R}^h: P_1(\underline{\alpha})\geq
0,\ldots,P_r(\underline{\alpha})\geq 0,
\ldots,P_s(\underline{\alpha})\geq 0\}$$
$$W_2=\{\underline{\alpha}\in{\R}^h: P_1(\underline{\alpha})\geq
0,\ldots,P_r(\underline{\alpha})\geq 0,
P_{r+1}(\underline{\alpha})\leq 0,\ldots,P_s(\underline{\alpha})\leq
0\}$$

Now we consider the polynomial:
$$R(\underline{x})=\sum_{i=r+1}^s P_i(\underline{x})$$ The
description of $W$ as union of $W_1$ and $W_2$ allows to conclude 
that inside $W$ an equation for $W_1$ is $R(\underline{x})\geq 0$ and
the equation for $W_2$ is $R(\underline{x})\leq 0$:
$$W_1=\{\underline{\alpha}\in W:R(\underline{\alpha})\geq 0\}
\qquad\quad W_2=\{\underline{\alpha}\in W:R(\underline{\alpha})\leq
0\}$$ what implies the following description for $W_1\cap W_2$:
$$W_1\cap W_2=\{\underline{\alpha}\in W:R(\underline{\alpha})=0\}$$
On $W_1\cap W_2$ we have $f=l_1=l_2$ and every zero of
$R(\underline{x})$ in
$W$ is a  zero of  $l_1(\underline{x})-l_2(\underline{x})$. So
applying \L ojasievicz  Inequality we obtain the existence of positive
integers  
$t$ and  $k$, and a positive number $c\in{\R}$ verifying:
$$\vert l_1(\underline{\alpha})-l_2(\underline{\alpha})\vert^t\leq
c\vert R(\underline{\alpha})\vert  (1+\Vert
\underline{\alpha}\Vert^2)^k\qquad\quad 
\forall \underline{\alpha}\in W$$ This allows to define the function:
$$v_{1,2}(\underline{\alpha})= l_2(\underline{\alpha})+\root t \of
{\max\{0,  c  R(\underline{\alpha})  (1+\Vert
\underline{\alpha}\Vert^2)^k\}}$$ verifying the desired conditions:

\smallskip \noindent [{$\bullet$}] for all $\underline{\alpha}\in W_2$ we have 
$v_{1,2}(\underline{\alpha})=l_2(\underline{\alpha})$,

\smallskip \noindent [{$\bullet$}] for all $\underline{\alpha}\in W_1$ we have:
$$v_{1,2}(\underline{\alpha})\geq l_2(\underline{\alpha})+
\vert l_1(\underline{\alpha})-l_2(\underline{\alpha})\vert
\geq l_1(\underline{\alpha})$$

Once all the functions $v_{n,p}$ have been constructed then it is 
very easy to check that:
$$f(\underline{\alpha})=\min\{\max\{v_{n,p}(\underline{\alpha}):
n\neq p, n\in N\}:p\in N\}\eqno(\star)$$ and the proof of the
theorem is obtained.\eop

\noindent Laureano Gonz\'alez--Vega ({\bf mailing address}),
Departamento de Matematicas, Estadistica y Computacion, Facultad de
Ciencias, Avenida de los Castros s/n, Universidad de Cantabria,
Santander 39071, Cantabria, Spain (e-mail: {\tt laureano.gonzalez@unican.es}).
 
 \medskip
\noindent Henri Lombardi, Laboratoire de Math\'ematiques, URA CNRS 741,
Universit\'e de Franche--Comt\'e, 25030 Besan{\c c}on, France (e-mail: {\tt henri.lombardi@univ-fcomte.fr}).
 
 \medskip
\noindent Louis Mah\'e, D\'epartement de Math\'ematiques, Universit\'e de
Rennes I, 35042 Rennes, France.

\end{document}